\documentclass[10pt]{amsart}
\usepackage{mathptmx}
\usepackage{amsmath}
\usepackage{amssymb}
\usepackage{array}
\usepackage{geometry}
\usepackage[bookmarks=true,colorlinks=true, pdfstartview=FitV, linkcolor=black, citecolor=blue, urlcolor=black]{hyperref}

\usepackage{color}
\definecolor{DarkRed}{rgb}{0.55,.00,0.2}
\definecolor{DarkGrey}{rgb}{0.35,.35,0.35}

\theoremstyle{definition}

\theoremstyle{remark}

\numberwithin{equation}{section}

\begin{document}

\title{On the generalized Dixon integral equation}

\author{Semyon Yakubovich$^\ast$}
\thanks{E-mail: syakubov@fc.up.pt}
\maketitle

\markboth{\rm \centerline{ Semyon   Yakubovich}}{}
\markright{\rm \centerline{Dixon equation  }}
\begin{center}{\it Department of Mathematics,  Faculty of Sciences,\\   University of Porto,  Campo Alegre st., 687,  4169-007 Porto,    Portugal}\end{center}

\begin{abstract} {\noindent A  unique analytic solution of the generalized Dixon nonhomogeneous integral equation is derived in $C^{1}[0,A],\ A > 0$.  It is written in terms of the Neumann series, which is expressed as a  double series of residues at  multiple poles of  powers of the gamma-function.}
\end{abstract}

\vspace{6pt}


\vspace{6pt}


\vspace{1cm} 

Consider the following nonhomogeneous  integral equation
$$f(x)= 1+  \lambda \int_0^A K(x,y) f(y) dy,\ x \in (0, A),\  \lambda \in \mathbb{C}\backslash \left\{0\right\},\  A \ge 1,\eqno(1)$$
where
$$K(x,y)= {\partial \over \partial y } P\left( {y\over x}\right)$$
and 
$$P(x)= {1\over B(a,a+1)}\int_0^x  {t^{a-1}\over (1+t)^{1+2a} }\ dt,\  x \ge 0,\ a >0,$$
where is the Euler beta-function \cite{erd}, Vol. I.  This equation is a generalization of the familiar Dixon integral equation \cite{tit},  Chapter 11.  Let $ f \in C^{1}[0, A]$.   Differentiating both sides of (1) with respect to $x$ and then integrating by parts, we find for $x >0$ 
$$f^\prime (x)=  \lambda {d\over dx} \left[   P\left( {y\over x}\right) f(y) \left|_0^A \right.  -  \int_0^A   P\left( {y\over x}\right)
 f^\prime (y) dy \right]$$
$$=  -  {\lambda \  x^{-a-1}   \over B(a,a+1)}  \left(1+  {y\over x}\right)^{-1-2a}   y^a f(y) \left|_0^A \right. $$
$$ + { \lambda\   x^{-a-1} \over  B(a,a+1) } \int_0^A    \left({x\over x+  y} \right)^{1+2a}   y^{a}  f^\prime (y) dy,\eqno(2) $$
where the differentiation under the integral sign is possible for all $x\ge 0$ due to the absolute and uniform convergence.  Indeed, we have 
$$ x^{-a-1} \int_0^A    \left({x\over x+  y} \right)^{1+2a}   y^{a}  | f^\prime (y)| dy \le \max_{x \in [0,A]} |f^\prime (x)|   \int_0^{A/x}   {y^{a} \over (1+y)^{1+2a}}  dy$$
$$ \le  \max_{x \in [0,A]} |f^\prime (x)|   \int_0^\infty   {y^{a} \over (1+y)^{1+2a}}  dy < \infty.\eqno(3)$$
Moreover, the integrated term vanishes in $y=0$.   Therefore we get from (2) the equation
$$ f^\prime (x) =   -  {\lambda f(A) A^a    \over B(a,a+1)}  {x^a\over \left(x+  A\right)^{1+2a} } $$
$$+   { \lambda\   \over  B(a,a+1) } \int_0^A  \left({x\over y} \right)^{a}   \left(1+ {x\over  y} \right)^{-1-2a}   f^\prime (y) {dy\over y}.\eqno(4) $$
Evidently, $f(x)\equiv 0$  does not satisfy the equation (1).  Moreover, we will show that $f(x) \equiv  \hbox{const}$ is not a solution of (1).  Indeed,  since  (1) can be written in the form 
$$ f(x)= 1  +  {\lambda   x^{a+1} \over B(a,a+1)} \int_0^A  {y^{a-1} \over   \left(x+  y\right)^{1+2a} }  f(y) dy, \eqno(5)$$
assuming that $f(x)= C$ is a solution, we have from (5)
$$C-1 =  {\lambda   C  \over B(a,a+1)} \int_0^{A/x}   {y^{a-1} \over   \left(1+  y\right)^{1+2a} }  dy,\ x >0.$$
Differentiating both sides of the latter equality with respect to $x$, we find 
$$-  {\lambda   C  \ A \over B(a,a+1) x^2}  \left({A\over x}\right)^{a-1}  \left(1+  {A\over x}\right)^{-1-2a}  =0$$
for all $x > 0$, which means that $C=0$.   But as we saw this is impossible.   Further, similar to (3) using the Weierstrass test of the uniform convergence, we verify that the integral in (5) converges absolutely and uniformly for all $x \ge 0$. Hence we find the limit values $f(0)= 1$ and $\lim_{x\to+\infty} f(x)= 1$.  

Let $g(x)=  f^\prime (x)  H(A-x),$ where $H(x)$ is the Heaviside function. Then applying the Mellin transform  \cite{tit}, \cite{yal},\  \cite{prud}, Vol. 3
$$g^*(s)= \int_0^\infty g(x) x^{s-1} dx,$$
to both sides of the equation (4), we take into account that the integral in (4) is the convolution for the Mellin transform of $g(x)$ \cite{tit} with the function 
$$h(x)=  {x^{a}\over (1+x)^{1+2a} }, $$
whose Mellin's transform is calculated via   the known beta-integral  \cite{prud}, Vol. I 
$$ \int_0^\infty   {x^{a-1}\over (1+x)^{1+2a} } \ dx = B(a, a+1),\ a > 0.\eqno(6)$$
Thus  we obtain $h^*(s)=  B(a+s, a+1-s).$ Consequently,   after application to (4) of  the Mellin transform and  the change of the order of integration by Fubini's  theorem by virtue of the estimate
$$ \int_0^\infty   x^{\sigma-1} dx  \int_0^\infty   \left({x\over y} \right)^{a}   \left(1+ {x\over  y} \right)^{-1-2a}  |g(y)| {dy\over y}$$
$$=   \int_0^\infty   {x^{\sigma+a -1} \over (1+x)^{1+2a} } dx  \int_0^A   y^{\sigma -1}  |f^\prime (y)| dy  $$ 
$$ \le    \max_{x \in [0,A]} |f^\prime (x)|  \int_0^\infty   {x^{\sigma+a -1} \over (1+x)^{1+2a} } dx  \int_0^A   y^{a -1}  dy 
 < \infty $$ 
where $\sigma = {\rm Re\   s} $,  $- a <  \sigma < a+1$, we come up with the algebraic  equality 
$$g^*(s) \left[ 1-    { \lambda\   B(a+s, a+1-s)  \over  B(a,a+1) }  \right] $$
$$ = -    \lambda f(A) A^{s-1} \  {   B(a+s, a+1-s)   \over B(a,a+1)} ,\quad  - a <  {\rm Re\   s}   < a+1.\eqno(7) $$
 Therefore,  if 
 $$\left|  { \lambda\   B(a+s, a+1-s)  \over  B(a,a+1) } \right| \le  {| \lambda| \   B(a+\sigma, a+1-\sigma)  \over  B(a,a+1) } < 1
 $$
i.e. 
 $$ |\lambda| <    { B(a,a+1)  \over B(a+\sigma, a+1-\sigma)  } \eqno(8)$$
 equality (7) becomes 
$$ g^*(s) =  -    \lambda f(A) A^{s-1} \  \frac {B(a+s, a+1-s)}{  B(a,a+1)  -  \lambda B(a+s, a+1-s) } .\eqno(9)$$
Now, considering $- a < \sigma <  a+1$, we use the inverse Mellin transform \cite{tit} to provide the unique solution of the equation (4) in the form
$$f^\prime (x)= - {  \lambda f(A)\over 2\pi A i}\int_{\sigma -i\infty}^{ \sigma + i\infty} 
 \frac {B(a+s, a+1-s)}{  B(a,a+1)  - \lambda  B(a+s, a+1-s) } \left({ x\over A} \right)^{-s} ds .\eqno(10)$$
Let $\sigma \neq 1$. Then  if $-a <\sigma < 1$, we integrate both sides of (10) from $0$ to $x$,  and taking into account account the value $f(0)=1$, we come up with the  unique soliton of the original integral equation (1), namely
 $$f (x)=  1 -  {  \lambda f(A) \over 2\pi i}\int_{\sigma -i\infty}^{ \sigma + i\infty} 
 \frac {B(a+s, a+1-s)}{  B(a,a+1)  -   \lambda B(a+s, a+1-s) } \left({ x\over A} \right)^{1 -s} {ds\over 1-s}.\eqno(11)$$

  Our final goal is to write the solution (11) in a different form, using the corresponding Neumann series and calculating the Mellin-Barnes integral by the residue theorem  at  multiple poles of powers of the  gamma-function \cite{erd}, Vol. 1,\ \cite{prud}, Vol. 3.  To do this, we will employ the series representation of the gamma-function $\Gamma(z)$ near the poles $z= -m,\  m \in \mathbb{N}_0$   (see, for instance, in \cite{eric} ) 
 $$\Gamma(z)=  {(-1)^m\over z+m} \sum_{k=0}^\infty c_{k,m}  (z+m)^k, \eqno(12)$$
 where 
$$ c_{k,m} = \sum_{\nu=0}^k \frac{(-1)^{(\nu+k)/2-1} (2^{k-\nu} - 2) B_{k-\nu} \pi^{k-\nu} }{\nu! (k-\nu)!} d_{\nu,m},\eqno(13)$$
$$d_{\nu, m} = \lim_{w\to m+1} {d^\nu\over d w^\nu} \left[ {1\over \Gamma(w) }\right] $$
and $B_\mu$ are the Bernoulli numbers \cite{erd}, Vol. I.  We begin, writing (11) in the form 
 $$f (x)= 1-  {  f(A)  \over 2\pi i}\int_{\sigma -i\infty}^{ \sigma + i\infty}  \sum_{n=1}^\infty
 \left({ \lambda  \   B(a+s, a+1-s) \over  B(a,a+1) } \right)^n   \left({ x\over A} \right)^{1 -s} {ds\over 1-s} .\eqno(14)$$
The change of the order of integration and summation is allowed owing to the absolute convergence which, in turn, can be verified  by virtue of the Stirling asymptotic formula for the gamma-function \cite{erd}, Vol. I , condition (7)  and the estimate 
$$   \int_{\sigma -i\infty}^{ \sigma + i\infty}   \sum_{n=1}^\infty \left|  \left({ \lambda  \   B(a+s, a+1-s) \over  B(a,a+1) } \right)^n   \left({ x\over A} \right)^{1 -s} {ds\over 1-s} \right| $$
$$\le     \left({ x\over A} \right)^{1 -\sigma}  \sum_{n=0}^\infty  \left({ |\lambda|  \    B(a+\sigma , a+1-\sigma) \over  B(a,a+1) } \right)^n    \int_{\sigma -i\infty}^{ \sigma + i\infty}  |  B(a+s, a+1-s) |    {|ds|\over |1-s|} < \infty, $$ 
 where $\sigma \in (-a, 1)$.  Therefore we write solution (14)  as 
 $$ f(x)=  1-  f(A)   \sum_{n=1}^\infty  \left({ \lambda  \over  \Gamma(a) \Gamma(a+1)} \right)^n{  1 \over 2\pi i}\int_{\sigma -i\infty}^{ \sigma + i\infty}      \left[ \Gamma (a+s) \Gamma(a+1-s)\right]^n  \left({ x\over A} \right)^{1 -s} {ds\over 1-s}.\eqno(15)$$
Meanwhile, the integral in (15) can be calculated by a special version of the residue theorem, namely, the Slater theorem (see \cite{prud}, Vol. 3) as the series of residues of the integrand at multiple left-hand poles $s= -a-m,\ m \in \mathbb{N}_{0}$ of the powers $[\Gamma (a+s)]^{n},\  n \in \mathbb{N}$.      We have
$${  1 \over 2\pi i}\int_{\sigma -i\infty}^{ \sigma + i\infty}      \left[ \Gamma (a+s) \Gamma(a+1-s)\right]^n   \left({ x\over A} \right)^{1 -s} {ds\over 1-s}  $$
$$= \sum_{m=0}^\infty  {\rm Res}_{s=-a-m} \left[\left[ \Gamma (a+s) \Gamma(a+1-s)\right]^n   \left({ x\over A} \right)^{1 -s}  {1\over 1-s} \right ].\eqno(16)$$
However,
$${\rm Res}_{s=-a-m} \left[\left[ \Gamma (a+s) \Gamma(a+1-s)\right]^n   \left({ x\over A} \right)^{1 -s} {1\over 1-s} \right ] $$
$$= {1\over   (n-1)!} \lim_{s\to -a-m} {d^{n-1} \over ds^{n-1} }\left[\left[ (s+a+m)\Gamma (a+s) \right]^n \left[ \Gamma(a+1-s)\right]^n   \left({ x\over A} \right)^{1 -s} {1\over 1-s}  \right ]  $$
$$=   {1\over  (n-1) !} \lim_{s\to -a-m} \sum_{r=0}^{n-1}  { n-1 \choose r } \  \left[\left[ (s+a+m)\Gamma (a+s) \right]^n  \right]^{(r)} \left[ \left[ \Gamma(a+1-s)\right]^n \left({ x\over A} \right)^{1 -s}  {1\over 1-s}    \right]^{(n-1-r)}.\eqno(17)$$
In the meantime,
$$\left[ \left[ \Gamma(a+1-s)\right]^n  \left({ x\over A} \right)^{1 -s} {1\over 1-s}   \right]^{(n-1-r)} $$
$$=  \sum_{\nu  =0}^{n-1-r}  { n-1-r \choose \nu }   \left( \left[ \Gamma(a+1-s)\right]^n \right)^{(\nu)} 
\left( \left({ x\over A} \right)^{1 -s} {1\over 1-s}  \right)^ {(n-1-r-\nu)}$$
$$=   \left({ x\over A} \right)^{1 -s} (n-1-r)! \sum_{\nu  =0}^{n-1-r} { (-1)^{n-1-r-\nu} \over \nu!}   \left( \left[ \Gamma(a+1-s)\right]^n \right)^{(\nu)} $$
$$\times \sum_{k=0}^{n-1-r-\nu}  {(-1)^k \left( \log (x/A) \right)^{n-1-r-\nu-k }  \over (n-1-r -\nu -k)!  (1-s)^{k+1} }.\eqno(18)$$
Meanwhile,  in order to calculate  higher -order derivatives  of  powers  $\left[ (s+a+m)\Gamma (a+s) \right]^n, \   \left[ \Gamma(a+1-s)\right]^n$  we will appeal to the familiar  Fa\' {a} di Bruno formula \cite{bruno}.   Thus we obtain 
$$\left[\left[ (s+a+m)\Gamma (a+s) \right]^n  \right]^{(r)}  =  \sum { r! \  n! \  \left( (s+a+m)\Gamma (a+s) \right)^{n-l}  \over  (n-l)!  b_1! b_2!\dots b_r! }  \left( {\left((s+a+m)\Gamma (a+s) \right)^{(1)} \over 1!} \right) ^{b_1} \dots$$
$$\times   \left( {\left((s+a+m)\Gamma (a+s)\right) ^{(r)}) \over r! }\right) ^{b_r},\eqno(19)$$ 
where the sum (19) is  over all different solutions in nonnegative integers $b_1,b_2,\dots, b_r$ of  
$b_1+2b_2+\dots+ rb_r =r$ and $ l = b_1+b_2+\dots + b_r$ and 
$$\left[\left[\Gamma (a+1- s) \right]^n  \right]^{(\nu)}  =  \sum { \nu ! \  n! \  \left(\Gamma (a+1-s) \right)^{n-l}  \over  (n-q)!  \beta_1! \beta_2!\dots \beta_\nu! }  \left( {\Gamma^{(1)}  (a+1- s)  \over 1!} \right) ^{\beta_1} \dots  \left( {\Gamma^{(\nu)}  (a+1-s) \over \nu! }\right) ^{b_\nu},\eqno(20)$$ 
where the sum  (20) is  over all different solutions in nonnegative integers $\beta_1, \beta_2,\dots, \beta_\nu$ of  
$\beta_1+2\beta_2+\dots+ r\beta_\nu  =\nu$ and $ q = \beta_1+\beta_2+\dots + \beta_\nu$.    Hence, recalling the expansion (12) of the gamma-function near the poles, we find, correspondingly, 
$$\lim_{s\to -a-m}  \left[\left[ (s+a+m)\Gamma (a+s) \right]^n  \right]^{(r)}  =  \sum { r! \  n! \  (-1)^{mn} \  c^{b_1}_{1,m}    \dots   c^{b_r}_{r,m} \over  (n-l)!  (m!)^{n-l} b_1! b_2!\dots b_r! },\eqno(21)$$ 
$$\lim_{s\to -a-m} \left[\left[\Gamma (a+1- s) \right]^n  \right]^{(\nu)}  =  \sum { \nu! \  n! \  \left(\Gamma (2a+1+m) \right)^{n-q}  \over  (n-q)!  \beta_1! \beta_2!\dots \beta_\nu ! }  \left( {\Gamma^{(1)}  (2a+1+m)  \over 1!} \right) ^{\beta_1} \dots  \left( {\Gamma^{(\nu)}  (2a+1+m) \over \nu! }\right) ^{\beta_\nu}.\eqno(22)$$ 
Consequently, combining with (17), (18), we obtain values of the residues 
$${\rm Res}_{s=-a-m} \left[\left[ \Gamma (a+s) \Gamma(a+1-s)\right]^n  \left({ x\over A} \right)^{1 -s} {1\over 1-s} \right] =  (n!)^2 \left({ x\over A} \right)^{a+1+m} \sum_{r=0}^{n-1}   \sum_{\nu  =0}^{n-1-r}  (-1)^{n-1-r-\nu}$$
$$\times     \sum_{k=0}^{n-1-r-\nu}  { (-1)^k \left( \log (x/A) \right)^{n-1-r-\nu-k }  \over (n-1-r -\nu -k)!  (a+1+m)^{k+1} }   \sum_{b} \sum_{\beta}  {  (-1)^{mn} \    c^{b_1}_{1,m}   \dots   c^{b_r}_{r,m}  \over  (n-l)!  (m!)^{n-l} b_1! b_2!\dots b_r! }  {   \Gamma^{n-q}  (2a+1+m)  \over  (n-q)!  \beta_1! \beta_2!\dots \beta_\nu ! }  $$
$$\times   \left( {\Gamma^{(1)}  (2a+1+m)  \over 1!} \right) ^{\beta_1} \dots  \left( {\Gamma^{(\nu)}  (2a+1+m) \over \nu! }\right) ^{\beta_\nu}.\eqno(23)$$
Thus the value of the integral (16) is equal to
$${  1 \over 2\pi i}\int_{\sigma -i\infty}^{ \sigma + i\infty}      \left[ \Gamma (a+s) \Gamma(a+1-s)\right]^n  \left({ x\over A} \right)^{1 -s} {ds\over 1-s}  $$
$$= (n!)^2 \sum_{m=0}^\infty \sum_{r=0}^{n-1}   \sum_{\nu  =0}^{n-1-r}  (-1)^{n-1-r-\nu} $$
$$\times   \left({ x\over A} \right)^{a+1+m}   \sum_{k=0}^{n-1-r-\nu}  {(-1)^k\ \left( \log (x/A) \right)^{n-1-r-\nu-k }  \over (n-1-r -\nu -k)!  (a+1+m)^{k+1} }   $$
$$\times  \sum_{b} \sum_{\beta}  {  (-1)^{mn} \  c^{b_1}_{1,m}  \dots    c^{b_r}_{r,m}    \over  (n-l)!  (m!)^{n-l} b_1! b_2!\dots b_r! }  {   \Gamma^{n-q}  (2a+1+m)  \over  (n-q)!  \beta_1! \beta_2!\dots \beta_\nu ! } $$
$$\times   \left( {\Gamma^{(1)}  (2a+1+m)  \over 1!} \right) ^{\beta_1} \dots  \left( {\Gamma^{(\nu)}  (2a+1+m) \over \nu! }\right) ^{\beta_\nu}\eqno(24)$$
and therefore our unique solution of the equation (1) takes the form  for $x \in [0,A]$
$$ f(x)= 1-  f(A)  \left({ x\over A} \right)^{a+1} \sum_{n=1}^\infty   \sum_{m=0}^\infty \left({ (-1)^m \lambda  \over  \Gamma(a) \Gamma(a+1)} \right)^n \     \left({ x\over A} \right)^{m}  \sum_{r=0}^{n-1}  \   \sum_{\nu  =0}^{n-1-r}  (-1)^{n-1-r-\nu}  $$
$$\times   \sum_{k=0}^{n-1-r-\nu}  {(-1)^k \  \left( \log (x/A) \right)^{n-1-r-\nu-k }  \over (n-1-r -\nu -k)!  (a+1+m)^{k+1} }  \sum_{b} \sum_{\beta} { n \choose l }  { n \choose q}  {  l! \  q! \    c^{b_1}_{1,m}  \dots    c^{b_r}_{r,m} \    \Gamma^{n-q}  (2a+1+m)   \over (m!)^{n-l} \  b_1! b_2!\dots \  b_r!  \  \beta_1! \beta_2!\dots \beta_\nu ! }  $$
$$\times  \left( {\Gamma^{(1)}  (2a+1+m)  \over 1!} \right) ^{\beta_1} \dots  \left( {\Gamma^{(\nu)}  (2a+1+m) \over \nu! }\right) ^{\beta_\nu}.\eqno(25)$$

Further, suppose that  $1 < \sigma < a+1$. Then recalling  (10), we integrate from $x$ to $A$ to get the equality
$$f(x) - f(A) =   {  \lambda f(A)\over 2\pi  i}\int_{\sigma -i\infty}^{ \sigma + i\infty} 
 \frac {B(a+s, a+1-s)}{  B(a,a+1)  - \lambda  B(a+s, a+1-s) } \left[1- \left({ x\over A} \right)^{1-s} \right] {ds\over 1-s} .\eqno(26)$$
Hence, in particular, passing to the limit in (26) when $x\to+\infty$, we remind the value 
$\lim_{x\to+\infty} f(x)= 1$ to derive the equality
$$ f(A)\left[ 1+  {  \lambda \over 2\pi  i}\int_{\sigma -i\infty}^{ \sigma + i\infty} 
 \frac {B(a+s, a+1-s)}{  B(a,a+1)  - \lambda  B(a+s, a+1-s) }  {ds\over 1-s} \right] = 1.\eqno(27)$$
Now, following the same scheme as above, we extend our solution for $x \ge A$. Indeed, making a simple substitution, we have (see (17), (18)) 
$${  1 \over 2\pi i}\int_{\sigma -i\infty}^{ \sigma + i\infty}      \left[ \Gamma (a+s) \Gamma(a+1-s)\right]^n  \left[1-  \left({ x\over A} \right)^{1 -s}\right] {ds\over 1-s}  $$
$$=  {  1 \over 2\pi i}\int_{1-\sigma -i\infty}^{1- \sigma + i\infty}      \left[ \Gamma (a+s) \Gamma(a+1-s)\right]^n  \left[1-  \left({ x\over A} \right)^{s} \right] {ds\over s}  $$
$$=  \sum_{m=0}^\infty  {\rm Res}_{s= -a-m} \left[\left[ \Gamma (a+s) \Gamma(a+1-s)\right]^n \left[  1- \left({ x\over A} \right)^{s} \right] {1\over s} \right ]$$
$$=   -  \sum_{m=0}^\infty   (-1)^{mn}   \sum_{r=0}^{n-1}  \   \sum_{\nu  =0}^{n-1-r}   \left[{1\over (a+m)^{n-r-\nu}} -  \left({ x\over A} \right)^{-a- m}    \sum_{k=0}^{n-1-r-\nu}  { \left( \log (x/A) \right)^{n-1-r-\nu-k }  \over (n-1-r -\nu -k)!  (a+m)^{k+1} }  \right] $$
$$\times  \sum_{b} \sum_{\beta} { n \choose l }  { n \choose q}  {  l! \  q! \    c^{b_1}_{1,m}  \dots    c^{b_r}_{r,m} \    \Gamma^{n-q}  (2a+1+m)   \over (m!)^{n-l} \  b_1! b_2!\dots \  b_r!  \  \beta_1! \beta_2!\dots \beta_\nu ! }  \left( {\Gamma^{(1)}  (2a+1+m)  \over 1!} \right) ^{\beta_1} \dots  \left( {\Gamma^{(\nu)}  (2a+1+m) \over \nu! }\right) ^{\beta_\nu}$$
and  (26) becomes $( x\ge A)$ 
$$f(x)= f(A) \left[ 1-   \sum_{n=1}^\infty  \sum_{m=0}^\infty \left({ \lambda  (-1)^m \over  \Gamma(a) \Gamma(a+1)} \right)^n   \  \sum_{r=0}^{n-1}  \   \sum_{\nu  =0}^{n-1-r}   \left[{1\over (a+m)^{n-r-\nu}} \right.\right.$$
$$\left. \left. -  \left({ x\over A} \right)^{-a- m} \     \sum_{k=0}^{n-1-r-\nu}  { \left( \log (x/A) \right)^{n-1-r-\nu-k }  \over (n-1-r -\nu -k)!  (a+m)^{k+1} }  \right] \right.$$
$$\times  \left. \sum_{b} \sum_{\beta} { n \choose l }  { n \choose q}  {  l! \  q! \    c^{b_1}_{1,m}  \dots    c^{b_r}_{r,m} \    \Gamma^{n-q}  (2a+1+m)   \over (m!)^{n-l} \  b_1! b_2!\dots \  b_r!  \  \beta_1! \beta_2!\dots \beta_\nu ! }  \left( {\Gamma^{(1)}  (2a+1+m)  \over 1!} \right) ^{\beta_1} \dots  \left( {\Gamma^{(\nu)}  (2a+1+m) \over \nu! }\right) ^{\beta_\nu}\right] .$$

We summarize our results by the following

{\bf Theorem.}  {\it Let $ f \in C^{1}[0, A],\ A \ge 1,\ \lambda \in  \mathbb{C}\backslash \left\{0\right\}$, satisfying the condition
$$|\lambda| <    { B(a,a+1)  \over B(a+\sigma, a+1-\sigma)  },\  a >0,$$
where $\sigma \in (- a, 1+a) \backslash\{1\}$.  Then the integral equation
$$f(x)= 1+  \lambda \int_0^A K(x,y) f(y) dy,$$
where
$$K(x,y)= {\partial \over \partial y } P\left( {y\over x}\right)$$
and 
$$P(x)= {1\over B(a,a+1)}\int_0^x  {t^{a-1}\over (1+t)^{1+2a} }\ dt$$
has the unique solution represented by the formula
$$ f(x)= 1-  f(A)  \left({ x\over A} \right)^{a+1} \sum_{n=1}^\infty   \sum_{m=0}^\infty \left({ (-1)^m \lambda  \over  \Gamma(a) \Gamma(a+1)} \right)^n \     \left({ x\over A} \right)^{m}  \sum_{r=0}^{n-1}  \   \sum_{\nu  =0}^{n-1-r}  (-1)^{n-1-r-\nu}  $$
$$\times   \sum_{k=0}^{n-1-r-\nu}  {(-1)^k \  \left( \log (x/A) \right)^{n-1-r-\nu-k }  \over (n-1-r -\nu -k)!  (a+1+m)^{k+1} }  \sum_{b} \sum_{\beta} { n \choose l }  { n \choose q}  {  l! \  q! \    c^{b_1}_{1,m}  \dots    c^{b_r}_{r,m} \    \Gamma^{n-q}  (2a+1+m)   \over (m!)^{n-l} \  b_1! b_2!\dots \  b_r!  \  \beta_1! \beta_2!\dots \beta_\nu ! }  $$
$$\times  \left( {\Gamma^{(1)}  (2a+1+m)  \over 1!} \right) ^{\beta_1} \dots  \left( {\Gamma^{(\nu)}  (2a+1+m) \over \nu! }\right) ^{\beta_\nu},$$
where  summations  over $b$ and $\beta$ are explained in $(20)$ and $c_{j,m}$ are defined by $(13)$. Moreover, it can be extended analytically for $x \ge A$ by the equality}

$$f(x)= f(A) \left[ 1-   \sum_{n=1}^\infty  \sum_{m=0}^\infty \left({ \lambda  (-1)^m \over  \Gamma(a) \Gamma(a+1)} \right)^n   \  \sum_{r=0}^{n-1}  \   \sum_{\nu  =0}^{n-1-r}   \left[{1\over (a+m)^{n-r-\nu}} \right.\right.$$
$$\left. \left. -  \left({ x\over A} \right)^{-a- m} \     \sum_{k=0}^{n-1-r-\nu}  { \left( \log (x/A) \right)^{n-1-r-\nu-k }  \over (n-1-r -\nu -k)!  (a+m)^{k+1} }  \right] \right.$$
$$\times  \left. \sum_{b} \sum_{\beta} { n \choose l }  { n \choose q}  {  l! \  q! \    c^{b_1}_{1,m}  \dots    c^{b_r}_{r,m} \    \Gamma^{n-q}  (2a+1+m)   \over (m!)^{n-l} \  b_1! b_2!\dots \  b_r!  \  \beta_1! \beta_2!\dots \beta_\nu ! }  \left( {\Gamma^{(1)}  (2a+1+m)  \over 1!} \right) ^{\beta_1} \dots  \left( {\Gamma^{(\nu)}  (2a+1+m) \over \nu! }\right) ^{\beta_\nu}\right] .$$

\bigskip
\centerline{{\bf Acknowledgments}}
\bigskip
I like to thank A. Polunchenko and I. Petrov   for bringing this integral equation and its  important applications in Statistics to my attention.

\bibliographystyle{amsplain}
{}

\end{document}